\documentclass[12pt,twoside,leqno,openany]{amsart}

\pdfoutput=1

\usepackage{amsbsy,amscd,amsfonts,amsmath,amssymb,amsthm,color,
fancybox,fancyhdr,footmisc,graphics,graphicx,ifthen,mathrsfs,
multicol,pdfpages,rotating,times,wasysym,pifont,blkarray}
\usepackage[dvipsnames,svgnames,x11names]{xcolor}

\usepackage[all]{xy}
\usepackage[utf8]{inputenc}
\usepackage[T1]{fontenc}
\usepackage{mathtools}
\newtagform{EngelLie}[\scriptstyle]{$}{$}
\makeatletter\newcommand{\leqnomode}{\tagsleft@true}
\newcommand{\reqnomode}{\tagsleft@false}\makeatother

\newtheorem{Theorem}[equation]{Theorem}

\theoremstyle{definition}

\newtheorem{Definition-Notation}[equation]{D\'efinition-Notation}

\newcommand{\C}{\mathbb{C}}

\newcommand{\R}{\mathbb{R}}

\newcommand{\NN}{\text{\sc n}}

\newcommand{\RR}{\text{\sc r}}

\newcommand{\XX}{\text{\sc x}}
\newcommand{\YY}{\text{\sc y}}
\newcommand{\ZZ}{\text{\sc z}}

\newcommand{\maux}{{\text{\usefont{T1}{qcs}{m}{sl}m}}}

\newcommand{\Waux}{{\text{\usefont{T1}{qcs}{m}{sl}W}}}

\definecolor{blue}{cmyk}{1.,1.,0.,0.63}
\definecolor{red}{cmyk}{0.,1.,1.,0.63}
\definecolor{green}{cmyk}{1.,0.,1.,0.63}
\definecolor{black}{cmyk}{1.,1.,1.,1.}

\makeatletter
\renewcommand{\@fnsymbol}[1]
{\ensuremath{\ifcase#1\or $*$\or $**$\or $***$\or $****$\or $*****$
\else\@ctrerr\fi}}
\makeatother

\newcommand{\HEAD}[2]{%
\pagestyle{fancy}
\fancyhead[RO]{\tiny\sf\thepage}
\fancyhead[CO]{{\tiny\sf #1}}
\fancyhead[LE]{\tiny\sf\thepage}
\fancyhead[CE]{{\tiny\sf #2}}
\fancyfoot{}}

\numberwithin{equation}{section}

\newcommand{\Section}[1]{
\renewcommand{\thesection}{\bf\arabic{section}}
\section{#1}
\renewcommand{\thesection}{\arabic{section}}}

\newcommand{\style}[1]{{\sf #1}}

\newcommand{\Aff}{\style{Aff}}

\renewcommand{\cos}{\style{cos}}

\renewcommand{\dim}{\style{dim}}

\renewcommand{\Im}{\style{Im}}

\renewcommand{\lim}{\style{lim}}

\renewcommand{\Re}{\style{Re}}

\renewcommand{\sin}{\style{sin}}

\newcommand{\Span}{\style{Span}}

\newcommand{\Hall}{\Hall}

\newcommand{\isqrt}{{\scriptstyle{\sqrt{-1}}}}

\newcommand{\vf}{\vfill\end{document}}

\makeatletter
\def\hlinethick#1{%
\noalign{\ifnum0=`}\fi\hrule \@height #1 \futurelet
\reserved@a\@xhline}
\makeatother

\newlength{\myskip}
\begin{document}

\setcounter{section}{0}

$\:$


\begin{center}

{\large\bf Homogeneous $\mathfrak{C}_{2,1}$ 
Models\footnotemark[1]}\medskip
\label{homogeneous-C21-models}

\medskip

\bigskip\bigskip

Wei-Guo {\sc Foo}\footnotemark[2],
Jo\"el~{\sc Merker}\footnotemark[3],
Pawe{\l}~{\sc Nurowski}\footnotemark[4],
The-Anh~{\sc Ta}\footnotemark[3]

\footnotetext[1]{\,
This research was supported
in part by the Polish National Science Centre (NCN) 
via the grant number 2018/29/B/ST1/02583,
and by the Norwegian Financial Mechanism
2014--2021 via the project registration number 2019/34/H/ST1/00636.}

\footnotetext[2]{\,\,
Institute of Mathematics, Academia Sinica, Taipei, Taiwan,
{\bf fooweiguo@gate.sinica.edu.tw}} 

\footnotetext[3]{\,\,
Département de Mathématiques d'Orsay,
CNRS, Université Paris-Saclay, 91405 Orsay Cedex,
France, {\bf joel.merker@universite-paris-saclay.fr}, 
{\bf tatheanhdtvt@gmail.com}}

\footnotetext[4]{\,\,Centrum Fizyki Teoretycznej,
Polska Akademia Nauk, Al. Lotnik\'ow 32/46, 02-668 Warszawa, Poland,
{\bf nurowski@cft.edu.pl}}

\end{center}\bigskip

\begin{center}
\begin{minipage}[t]{12.5cm}
\parindent 0.53cm
\footnotesize
\noindent
{\sc Abstract}. 
Fels-Kaup (Acta Mathematica 2008) classified homogeneous
$\mathfrak{C}_{2,1}$ hypersurfaces $M^5 \subset \C^3$, and discovered
that they are {\em all} biholomorphic to {\em tubes} $S^2 \times
i\R^3$ over some affinely homogeneous surface $S^2 \subset \R^3$.
The second and third authors 2003.08166, by performing {\em highly
non-straightforward calculations}, conducted the {\sl Cartan method
of equivalence} to classify homogeneous models of PDE systems
related to such $\mathfrak{C}^{2,1}$ hypersurfaces $M^5 \subset \C^3$.

Kolar-Kossovskiy 1905.05629 and the authors 2003.01952
constructed a formal and a convergent 
Poincar\'e-Moser normal form for
$\mathfrak{C}_{2,1}$ hypersurfaces $M^5 \subset \C^3$.
But this was only a first, preliminary step.
Indeed, the invariant branching tree underlying
Fels-Kaup's classification was still missing in the literature, due to
computational obstacles.

The present work applies the {\sl power series method of
equivalence}, confirms Fels-Kaup 2008, and finds the 
{\em differential-invariant} tree:
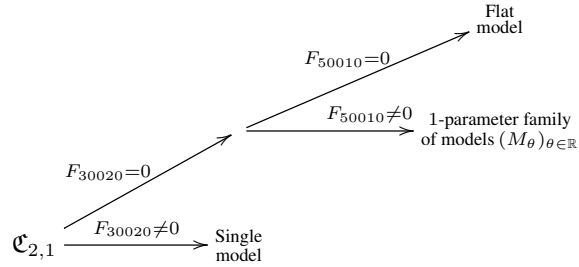
\[
\xymatrix{
&&
&&
{\substack{
\text{Flat}
\\
\text{model}}}
\\
&&
\ar[urr]^{F_{50010}=0}
\ar[rr]^{F_{50010}\neq0}
&&
{\substack{
1\text{-parameter family}
\\
\text{of models}\,(M_\theta)_{\theta\in\R}}}
\\
\ar[urr]^{F_{30020}=0}
\mathfrak{C}_{2,1}
\ar[rr]^{F_{30020}\neq0}
&&
{\substack{
\text{Single}
\\
\text{model}}}
&&
}
\]

To terminate the middle (thickest) branch, it is necessary
to compute up to order $10$ with $5$ variables.
Again, calculations, done by hand, are {\em non-straightforward}.
\end{minipage}
\end{center}

\Section{\bf Introduction}
\label{introduction-C21}
\HEAD{{\ref{introduction-C21}}.~{\sf Introduction}
}{
Wei-Guo {\sc Foo}, Jo\"el {\sc Merker}, Pawe{\l}~{\sc Nurowski},
The-Anh {\sc Ta}}

Strong interactions between affine and CR 
{\em homogeneous} geometries are know to exist
{\cite{Nurowski-Tafel-1988,
Nurowski-Sparling-2003,
Gaussier-Merker-2003,
Fels-Kaup-2007, 
Fels-Kaup-2008,
Isaev-2011,
Merker-Pocchiola-2018,
Isaev-2016,
Isaev-2016-bis,
Isaev-2018,
Foo-Merker-2019,
Merker-Nurowski-2019,
Loboda-2020,
Doubrov-Merker-The-2020,
Merker-Nurowski-2021}}.
Notably, {\em a CR homogeneous
$M \subset \C^\NN$ is biholomorphic to a tube
over an affinely homogeneous hypersurface $H \subset \R^\NN$ 
if and only if its infinitesimal Lie symmetry algebra 
contains a (nontrivial) maximally real Abelian ideal}.

Indeed, let $M \subset \C^{\NN \geqslant 2}$ be a local $\mathcal{C}^\omega$
CR hypersurface, in coordinates 
$\ZZ = (\ZZ_1, \dots, \ZZ_\NN) \in \C^\NN$, with $0 \in M$.
Assume that $M$ is CR-homogeneous, so that the {\em real} Lie algebra:
\[
\mathfrak{hol}(M) 
:= 
\Big\{ 
L
= 
{\textstyle{\sum_{i=1}^\NN}}\, 
a_i(\ZZ)\, 
\tfrac{\partial}{\partial\ZZ_i}\,
\text{holomorphic}
\colon\,\,
\big(L+\overline{L}\big)
\big\vert_M\,\,
\text{is tangent to}\,\,
M
\Big\},
\]
is of dimension $\dim\,M \leqslant \RR \leqslant \infty$,
due to $T_0 M = \Span \big\{ 
(L+\overline{L}) \big\vert_0 \colon 
L \in \mathfrak{hol}(M) \big\}$.

If $\mathfrak{hol}(M) \supset \mathfrak{a}$ contains
an $\NN$-dimensional Abelian (real) Lie subalgebra $\mathfrak{a} = 
\Span\, \big(L_1, \dots, L_\NN \big)$
of holomorphic vector fields having 
{\em maximally real} span: 
\[
\Span\,
\big(
L_1+\overline{L}_1\big\vert_0, 
\dots, 
L_\NN+\overline{L}_\NN 
\big\vert_0\big)
\,\subset\,
T_0\C^\NN,
\]
then after a straightening biholomorphism, one has
$L_1 = \isqrt\, \partial_{\ZZ_1}$, \dots, $L_\NN = 
\isqrt\, \partial_{\ZZ_\NN}$. 

Assume furthermore that
$\mathfrak{a} \subset \mathfrak{hol}(M)$ is an {\em ideal}.
Consider other $L_\nu \in \mathfrak{hol}(M)$ 
for $\NN+1 \leqslant \nu \leqslant \RR$ completing a basis. 
Since each $\big[ \isqrt\, \partial_{\ZZ_i}, L_\nu \big]$
must be a real linear combination of $\isqrt\,\partial_{\ZZ_1}, 
\dots, \isqrt\, \partial_{\ZZ_\NN}$,
it comes: 
\[
L_\nu
\,=\,
\sum_{i=1}^\NN\,
\Big(
\sum_{j=1}^\NN\,
a_{\nu,i,j}\,\ZZ_j
+
b_{\nu,i}
\Big)\,
\frac{\partial}{\partial\ZZ_i}
\eqno
{\scriptstyle{(\NN+1\,\leqslant\,\nu\,\leqslant\,\RR)}},
\]
with constants $a_{\nu,i,j} \in \R$, and $b_{\nu, i} \in \C$;
in fact $b_{\nu, i} \in \R$,
after subtracting appropriate 
linear combinations of the $\isqrt\,\partial_{\ZZ_i}$.
Tangency to $M$ of the real parts of the $\isqrt\, \partial_{\ZZ_i}$
implies that $M = H \times i\R^\NN$ 
with $H \subset \R^\NN$ a hypersurface.
Furthermore, writing 
$\ZZ_i = \XX_i + \isqrt\, \YY_i$, the vector fields
\[
T_\nu
\,:=\,
\sum_{i=1}^\NN\,
\Big(
\sum_{j=1}^\NN\,
a_{\nu,i,j}\,\XX_j
+
b_{\nu,i}
\Big)\,
\frac{\partial}{\partial\XX_i}
\eqno
{\scriptstyle{(\NN+1\,\leqslant\,\nu\,\leqslant\,\RR)}},
\]
are tangent to $H$, and their span at $0 \in H$ spans $T_0 H$.
The converse is direct.\hfill$\bigtriangleup$

\smallskip

In addition, $\mathfrak{C}_{2,1}$\,\,---\,\,abbreviation for 
"{\sl $2$-nondegenerate of constant Levi rank 
$1$}"\,\,---\,\,hypersurfaces $M^5 \subset \C^3$
correspond 
{\cite{Fels-Kaup-2007, Fels-Kaup-2008,
Chen-Merker-2019}}
to parabolic surfaces $S^2 \subset \R^3$. 
The classification of all affinely 
homogeneous surfaces 
$S^2 \subset \R^3$ appears in~{\cite{Abdalla-Dillen-Vrancken-1997,
Doubrov-Komrakov-Rabinovich-1996, Eastwood-Ezhov-1999}}.
The noncylindrical parabolic ones can be presented as 
follows~{\cite{Fels-Kaup-2008}}:

\smallskip\noindent{\bf (1)}\,
$\big\{ x_1^2 + x_2^2 = x_3^2, \,\, 
x_3 >0 \big\}$;

\smallskip\noindent{\bf (2a)}\,
$\big\{ r ( \cos\,t, \sin\,t, e^{ \omega t})\in \mathbb{ R}^3: \, r
\in \mathbb{ R}^+ \ \text{\rm and}\ t \in \mathbb{ R} \big\}$ with
$\omega > 0$ arbitrary;

\smallskip\noindent{\bf (2b)}\,
$\big\{ r ( 1, t, e^t) \in \mathbb{ R}^3: \, r \in \mathbb{ R}^+
\ \text{\rm and}\ t \in \mathbb{ R} \big\}$;

\smallskip\noindent{\bf (2c)}\,
$\big\{ r ( 1, e^t, e^{ \theta t}) \in \mathbb{ R}^3 : \, r \in
\mathbb{ R}^+ \ \text{\rm and}\ t \in \mathbb{ R}
\big\}$ with $\theta > 2$
arbitrary;

\smallskip\noindent{\bf (3)}\,
$\big\{ c ( t) + r c' ( t) \in \mathbb{ R}^3 : \, r \in \mathbb{
R}^+ \ \text{\rm and} \ t \in \mathbb{ R} \big\}$,
where $c(t) = (t, t^2, t^3)$.

\smallskip

Fels-Kaup's striking result was: {\em Every locally homogeneous
$\mathfrak{C}_{2,1}$ hypersurface $M^5 \subset \C^3$ is locally
biholomorphic to $S^2 \times i\R^3$, with $S^2 \subset \R^3$ being
one of {\bf (1)}, {\bf (2a)}, {\bf (2b)}, {\bf (2c)}, {\bf (3)};
distinct such $S^2 \times i\R^3$ are pairwise biholomorphically
inequivalent; all but the (flat model) tube {\bf (1)} are simply
transitive}.

\smallskip

Fels-Kaup's proof relied on expert knowledge of Lie structure theory.
But only the equivalence method can reach information about CR
invariants. The goal of this article is to explore the concerned
CR invariants (either relative or absolute), 
since nothing about the branchings they create
appears in~{\cite{Fels-Kaup-2008, 
Medori-Spiro-2014,
Foo-Merker-2019,
Kolar-Kossovskiy-2019, 
Foo-Merker-Ta-2020}}.

\Section{\bf Affine Classification}
\label{affine-classification}
\HEAD{{\ref{affine-classification}}.~{\sf Affine Classification}
}{
Wei-Guo {\sc Foo}, Jo\"el {\sc Merker}, Pawe{\l}~{\sc Nurowski},
The-Anh {\sc Ta}}

Before stating CR results, 
let us present an alternative
(elementary) classification of $\Aff(\R^3)$-homogeneous noncylindrical
parabolic surfaces, whose final invariant tree is:
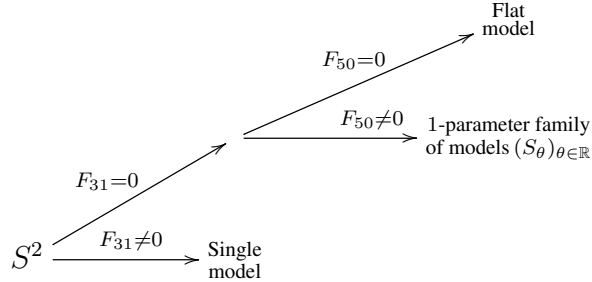
\[
\xymatrix{
&&
&&
{\substack{
\text{Flat}
\\
\text{model}}}
\\
&&
\ar[urr]^{F_{50}=0}
\ar[rr]^{F_{50}\neq0}
&&
{\substack{
1\text{-parameter family}
\\
\text{of models}\,(S_\theta)_{\theta\in\R}}}
\\
\ar[urr]^{F_{31}=0}
S^2
\ar[rr]^{F_{31}\neq0}
&&
{\substack{
\text{Single}
\\
\text{model}}}
&&
}
\]

To explain this tree, let $S^2 \subset \R^3$ with $0 \in S^2$ be
$\mathcal{C}^\omega$ graphed as $u = F(x,y) = \sum_{j+k \geqslant 1}\,
F_{j,k}\, \frac{1}{j!}\, x^j \frac{1}{k!}\, y^k$. 
According to~{\cite{Chen-Merker-2019}},
parabolicity expresses as $F_{xx} \neq 0 \equiv \big\vert
\begin{smallmatrix} F_{xx} & F_{xy} \\ F_{yx} & F_{yy} 
\end{smallmatrix} \big\vert$,
and non-cylindricality as $0 \neq F_{xxy}\, F_{xx} - F_{xxx} \,
F_{xy}$.

A preliminary normalization is:
\[
u
\,=\,
\tfrac{1}{2}\,
x^2
+
\tfrac{1}{2}\,
x^2y
+
\tfrac{1}{6}\,
F_{3,1}\,x^3y
+
\tfrac{1}{2}\,
x^2y^2
+
{\rm O}_{x,y}(5).
\]
The coefficient $F_{3,1}$ is a relative invariant under
$\Aff(\R^3)$, hence it creates a branching.

In the branch $F_{3,1} \neq 0$, one normalizes $F_{3, 1} := 1$,
and also $F_{4,1} := 0$, 
whence up to order $5$:
\[
\aligned
u
&
\,=\,
\tfrac{1}{2}\,
x^2
+
\tfrac{1}{2}\,
x^2y
+
\tfrac{1}{6}\,
x^3y
+
\tfrac{1}{2}\,
x^2y^2
\\
&
\ \ \ \ \
+
\tfrac{1}{120}\,F_{5,0}\,
x^5
+
\tfrac{1}{2}\,
x^3y^2
+
\tfrac{1}{2}\,
x^2y^3
+
{\rm O}_{x,y}(6).
\endaligned
\]
One finds that $F_{5,0} = \frac{20}{9}$ necessarily,
and that all higher order $F_{j,k}$ are uniquely 
determined {\em constants}, for instance up to order $8$:
\[
\aligned
u
&
\,=\,
\tfrac{1}{2}\,
x^2
+
\tfrac{1}{2}\,
x^2y
+
\tfrac{1}{6}\,
x^3y
+
\tfrac{1}{2}\,
x^2y^2
\\
&
\ \ \ \ \
+
\tfrac{1}{54}\,
x^5
+
\tfrac{1}{2}\,
x^3y^2
+
\tfrac{1}{2}\,
x^2y^3
\\
&
\ \ \ \ \
+
\tfrac{1}{162}\,
x^6
+
\tfrac{1}{18}\,
x^5y
+
\tfrac{1}{8}\,
x^4y^2
+
x^3y^3
+
\tfrac{1}{2}\,x^2y^4
\\
&
\ \ \ \ \
-\,
\tfrac{1}{486}\,x^7
+
\tfrac{7}{108}\,x^6y
+
\tfrac{5}{54}\,x^5y^2
+
\tfrac{5}{8}\,
x^4y^3
+
\tfrac{5}{3}\,
x^3y^4
+
\tfrac{1}{2}\,
x^2y^5
\\
&
\ \ \ \ \
+
\tfrac{5}{5832}\,x^8
+
\tfrac{1}{162}\,x^7y
+
\tfrac{1}{4}\,x^6y^2
+
\tfrac{47}{216}\,
x^5y^3
+
\tfrac{15}{8}\,
x^4y^4
+
\tfrac{5}{2}\,
x^3y^5
+
\tfrac{1}{2}\,
x^2y^6
+
{\rm O}_{x,y}(9).
\endaligned
\]
The (transitive) affine Lie symmetry algebra is $2$-dimensional,
generated by:
\[
\aligned
e_1
&
\,:=\,
\big(
1+x-y-\tfrac{10}{9}\,u
\big)\,
\partial_x
+
\big(
\tfrac{10}{9}\,x-y-\tfrac{10}{9}\,u
\big)\,
\partial_y
+
\big(
x+2\,u
\big)\,
\partial_u,
\\
e_2
&
\,:=\,
\big(
-2\,x+u
\big)\,
\partial_x
+
\big(
1-\tfrac{4}{3}\,x-y+\tfrac{8}{9}\,u
\big)\,
\partial_y
-
3\,u\,
\partial_u,
\endaligned
\]
having Lie bracket:
\[
[e_1,e_2]
\,=\,
-\,e_1
-
\tfrac{1}{3}\,e_2.
\]
This is {\bf (3)}.

Next, consider the (invariant) branch $F_{3,1} = 0$.
Necessarily, $F_{4,1} = 0$, hence:
\[
u
\,=\,
\tfrac{1}{2}\,x^2
+
\tfrac{1}{2}\,x^2y
+
\tfrac{1}{2}\,x^2y^2
+
\tfrac{1}{120}\,F_{5,0}\,
x^5
+
\tfrac{1}{2}\,
x^2y^3
+
{\rm O}_{x,y}(6),
\]
with $F_{5,0}$ being a relative invariant, again
creating a (sub)branching.

The (sub)branch $F_{3,1} = 0 = F_{5,0}$ 
conducts~{\cite{Chen-Merker-2019, Merker-Nurowski-2021}}
to the {\sl flat model}\,\,---\,\,a graphed representation
of {\bf (1)} above\,\,---:
\[
u
\,=\,
\frac{1}{2}\,
\frac{x^2}{1-y},
\]
having $4$-dimensional (transitive) affine
Lie symmetry algebra generated by:
\[
e_1
\,:=\,
(1-y)\,\partial_x
+
x\,\partial_u,
\ \ \ \ \ \ \
e_2
\,:=\,
(1-y)\,\partial_y
+
u\,\partial_u,
\ \ \ \ \ \ \
e_3
\,:=\,
x\,\partial_x
+
2\,u\,\partial_u,
\ \ \ \ \ \ \
e_4
\,:=\,
-\,u\,\partial_x
+
x\,\partial_y,
\]
with structure:
\[
[e_1,e_2]
\,=\,
e_1,
\ \ \ \ \
[e_1,e_3]
\,=\,
e_1,
\ \ \ \ \
[e_1,e_4]
\,=\,
e_2,
\ \ \ \ \
[e_2,e_4]
\,=\,
e_4,
\ \ \ \ \
[e_3,e_4]
\,=\,
e_4.
\]

In the thickest (sub)branch $F_{3,1} = 0 \neq F_{5,0}$, 
one normalizes $F_{5,0} := 1$, and also $F_{6,0} := 0$.
Necessarily, $F_{5,1} = 4$ and $F_{6,1} = 0$, hence:
\[
\aligned
u
&
\,=\,
\tfrac{1}{2}\,x^2
+
\tfrac{1}{2}\,x^2y
+
\tfrac{1}{2}\,x^2y^2
+
\tfrac{1}{120}\,
x^5
+
\tfrac{1}{2}\,
x^2y^3
\\
&
\ \ \ \ \
+
\tfrac{1}{30}\,
x^5y
+
\tfrac{1}{2}\,
x^2y^4
+
\tfrac{1}{5040}\,
F_{7,0}\,
x^7
+
\tfrac{1}{12}\,
x^5y^2
+
\tfrac{1}{2}\,
x^2y^5
+
{\rm O}_{x,y}(8),
\endaligned
\]
with $F_{7,0}$ being an {\em absolute} invariant.
Call it:
\[
F_{7,0}
\,=:\,
\theta.
\]

One therefore finds a $1$-parameter family of affinely
inequivalent homogeneous models $\big( S_\theta \big)_{\theta 
\in \R}$, with $2$-dimensional 
(simply transitive) affine Lie symmetry algebra:
\[
\aligned
e_1
&
\,:=\,
\big(
1-y+\tfrac{1}{3}\,\theta\,u
\big)\,
\partial_x
+
\big(
-\,\tfrac{1}{3}\,\theta\,x
-\tfrac{1}{6}\,u
\big)\,
\partial_y
+
x\,\partial_u,
\\
e_2
&
\,:=\,
-\,x\,\partial_x
+
(1-y)\,\partial_y
-
u\,\partial_u,
\ \ \ \ \ \ \ \ \ \ \ \ \ \ \ \ \ \ \ \ \ \ \ \ \ \ \ \ \ \ \ \ \ \ \ 
\ \ \ \ \ \ \ \ \ \ \ \ \ \ \ \ \
[e_1,e_2]
\,=\,
0.
\endaligned
\]
This unifies {\bf (2a)}, {\bf (2b)}, {\bf (2c)}.

\Section{\bf CR Classification}
\label{CR-classification}
\HEAD{{\ref{CR-classification}}.~{\sf CR Classification}
}{
Wei-Guo {\sc Foo}, Jo\"el {\sc Merker}, Pawe{\l}~{\sc Nurowski},
The-Anh {\sc Ta}}

In coordinates $\C^3 \ni \big(z, \zeta, w = u + \isqrt\, v\big)$, 
the graphed representation~{\cite{Gaussier-Merker-2003,
Fels-Kaup-2007,
Foo-Merker-Ta-2019,
Foo-Merker-Ta-2020,
Chen-Foo-Merker-Ta-2020}} 
of the flat model is due to Gaussier-Merker:
\[
u
\,=\,
\frac{z\overline{z}
+\frac{1}{2}\,\overline{z}^2\zeta
+\frac{1}{2}\,z^2\overline{\zeta}}{1-\zeta\overline{\zeta}}
\,\,=:\,\,
\maux\big(z,\zeta,\overline{z},\overline{\zeta}\big).
\]
The $5$-dimensional
Lie group of its automorphisms fixing the origin writes:
\[
\aligned
z'
&
\,:=\,
\lambda\,
\frac{z+i\,\alpha\,z^2+
\big(i\,\alpha\,\zeta-i\,\overline{\alpha}\big)\,w}
{1+2i\,\alpha\,z-\alpha^2z^2
-
\big(\alpha^2\zeta-\alpha\overline{\alpha}+i\,\rho\big)\,w},
\\
\zeta'
&
\,:=\,
\frac{\lambda}{\overline{\lambda}}\,\,
\frac{\zeta+2i\,\overline{\alpha}\,z
-\big(\alpha\overline{\alpha}+i\,\rho\big)\,z^2
+\big(\overline{\alpha}^2-i\,\rho\,\zeta
-\alpha\overline{\alpha}\,\zeta\big)\,w}
{1+2i\,\alpha\,z-\alpha^2z^2
-
\big(\alpha^2\zeta-\alpha\overline{\alpha}+i\,\rho\big)\,w},
\\
w'
&
\,:=\,
\lambda\overline{\lambda}\,
\frac{w}
{1+2i\,\alpha\,z-\alpha^2z^2
-
\big(\alpha^2\zeta-\alpha\overline{\alpha}+i\,\rho\big)\,w},
\endaligned
\]
where $\lambda \in \C^\ast$, $\alpha \in \C$, $\rho \in \R$ are
free.

A general $\mathfrak{C}_{2,1}$ hypersurface $M^5 \subset \C^3$
with $0 \in M$ writes as a perturbation of this model:
\[
u
\,=\,
F\big(z,\zeta,\overline{z},\overline{\zeta},v\big)
\,=\,
\maux(z,\zeta,\overline{z},\overline{\zeta})
+
G\big(z,\zeta,\overline{z},\overline{\zeta},v\big),
\]
where:
\[
F
\,=\,
\sum_{h,i,j,k,l}\,
z^h\zeta^i\overline{z}^j\overline{\zeta}^kv^l\,
F_{h,i,j,k,l}
\,=\,
\sum_{h,i,j,k}\,
z^h\zeta^i\overline{z}^j\overline{\zeta}^k
F_{h,i,j,k}(v),
\]
with $\overline{F_{h,i,j,k,l}} = F_{j,k,h,i,l}$,
with $0 = F_{0,0,0,0,0}$, 
and the same for $G$.

The Poincar\'e-Moser {\em convergent}
normal form established in~{\cite{Foo-Merker-Ta-2020}}
shows that, after some local biholomorphism fixing the origin,
one can assume:
\[
\aligned
0
&
\,\equiv\,
F_{h,i,0,0}(v),
\\
0
&
\,\equiv\,
F_{h,i,1,0}(v),
\\
0
&
\,\equiv\,
F_{h,i,2,0}(v),
\endaligned
\ \ \ \ \ \ \ \ \ \ \ \ \ \ \ \ \ \ \ \
\aligned
0
&
\,\equiv\,
F_{3,0,0,1}(v),
\\
0
&
\,\equiv\,
F_{4,0,0,1}(v)
\,\equiv\,
F_{3,0,1,1}(v),
\\
0
&
\,\equiv\,
F_{4,0,1,1}(v)
\,\equiv\,
F_{3,0,3,0}(v),
\endaligned
\]
with the exceptions $1 \equiv F_{1,0,1,0}(v)$ and 
$\frac{1}{2} \equiv F_{2,0,0,1}(v)$.

Suppose ${M'}^5 \subset {\C'}^3$ is another such $\mathfrak{C}_{2,1}$
hypersurface, similarly normalized. 
If: 
\[
(z,\zeta,w)
\,\longmapsto\,
\big(
f(z,\zeta,w),\,
g(z,\zeta,w),\,
h(z,\zeta,w)
\big)
\,=:\,
(z',\zeta',w'),
\]
is a local holomorphic map fixing the origin which sends $M$ 
into $M'$, then as follows from general Poincar\'e-Moser theory,
it is of the form above for certain five real parameters
$\lambda \in \C^\ast$, $\alpha \in \C$, $\rho \in \R$.
Our goal is to normalize this remaining ambiguity,
{\em cf.} Questions~{\bf Q}\textsuperscript{\ding{192}}
and~{\bf Q}\textsuperscript{\ding{195}} 
in~{\cite{Foo-Merker-Ta-2020}}.

Attributing weights $[z] := 1$, $[\zeta] := 1$,
$[w] := 2$, let us therefore show weighted order $5$ terms:
\[
\aligned
u
&
\,=\,
z\overline{z}
+
\tfrac{1}{2}\,
\overline{z}^2\zeta
+
\tfrac{1}{2}\,
z^2\overline{\zeta}
+
z\overline{z}\zeta\overline{\zeta}
+
\tfrac{1}{2}\,
\overline{z}^2\zeta\zeta\overline{\zeta}
+
\tfrac{1}{2}\,
z^2\overline{\zeta}\zeta\overline{\zeta}
\\
&
\ \ \ \ \
+
2\,\Re\,
\Big\{
z^3\overline{\zeta}^2\,
F_{3,0,0,2,0}
\Big\}
+
{\rm O}_{z,\zeta,\overline{z},\overline{\zeta},v}(6),
\endaligned
\]
the remainder being {\em weighted} as well.
This coefficient $F_{3,0,0,2,0}$ is a relative invariant,
hence it creates a branching.
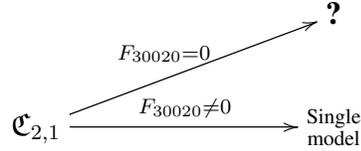
\[
\xymatrix{
&&&
{\text{\bf ?}}
\\
\ar[urrr]^{F_{30020}=0}
\mathfrak{C}_{2,1}
\ar[rrr]^{F_{30020}\neq0}
&&&
{\substack{
\text{Single}
\\
\text{model}}}
}
\]

\begin{Theorem}
In the branch $F_{3,0,0,2,0} \neq 0$, one can normalize
$F_{3,0,0,2,0} := 1$, so $\lambda := 1$, and
$3$ supplementary (real) normalizations hold:
\[
\aligned
F_{4,0,0,2,0}
&
\,:=\,
0,
&
\ \ \ \ \ \ \ \ \ \ \ \ \ \ \ \ \ \ \ \
\text{so}\ \
\alpha
&
\,:=\,
0,
\\
\Im\,F_{3,0,2,1,0}
&
\,:=\,
0,
&
\ \ \ \ \ \ \ \ \ \ \ \ \ \ \ \ \ \ \ \
\text{so}\ \
\rho
&
\,:=\,
0,
\endaligned
\]
so that the isotropy is reduced to be zero-dimensional.

Furthermore, all coefficients $F_{h,i,j,k,l} \in \C$ are
uniquely determined to be specific constants, with:
\[
F
\,=\,
F^2+F^3+F^4+F^5+F^6+F^7+F^8
+
{\rm O}_{z,\zeta,\overline{z},\overline{\zeta},v}(9),
\]
where:
\[
\aligned
F^2
&
\,=\,
z\overline{z},
\\
F^3
&
\,=\,
\tfrac{1}{2}\,
\overline{z}^2\zeta
+
\tfrac{1}{2}\,
z^2\overline{\zeta},
\\
F^4
&
\,=\,
z\overline{z}\zeta\overline{\zeta},
\\
F^5
&
\,=\,
z^3\overline{\zeta}^2
+
\zeta^2\overline{z}^3
+
\tfrac{1}{2}\,
z^2\zeta\overline{\zeta}^2
+
\tfrac{1}{2}\,
\zeta^2\overline{z}^2\overline{\zeta},
\endaligned
\]
\[
\aligned
F^6
&
\,=\,
-\,
2\,z\zeta^2\overline{z}^3
-
2\,z^2\zeta\overline{z}^3
+
3\,z^2\zeta\overline{z}\overline{\zeta}^2
+
\tfrac{1}{3}\,
\zeta^3\overline{z}^3
+
z\zeta^2\overline{z}\overline{\zeta}^2
-
2\,z^3\overline{z}^2\overline{\zeta}
+
3\,z\zeta^2\overline{z}^2\overline{\zeta}
\\
&
\ \ \ \ 
+\tfrac{1}{3}\,z^3\,\overline{\zeta}^3
-
2\,z^3\overline{z}\overline{\zeta}^2,
\endaligned
\]
\[
\aligned
F^7
&
\,=\,
\tfrac{8}{3}\,\zeta^3\overline{z}^4
-
6\,z^2\zeta^2\overline{z}^2\overline{\zeta}
+
z\zeta^3\overline{z}^2\overline{\zeta}
+
4\,z^3\zeta\overline{z}^3
+
\tfrac{8}{3}\,z^4\overline{\zeta}^3
+
z^4\overline{z}^3
+
z^3\overline{z}^4
\\
&
\ \ \ \ \
+
3\,\zeta^3\overline{z}^3\overline{\zeta}
-
\tfrac{4}{3}\,z^3\overline{z}\overline{\zeta}^3
-
2\,z^4\overline{z}\overline{\zeta}^2
+
z^2\zeta\overline{z}\overline{\zeta}^3
-
\tfrac{4}{3}\,z\zeta^3\overline{z}^3
+
2\,z^3\overline{z}^2\overline{\zeta}^2
\\
&
\ \ \ \ \
+
\tfrac{1}{2}\,\zeta^3\overline{z}^2\overline{\zeta}^2
+
3\,z^3\zeta\overline{\zeta}^3
-
6\,z^2\zeta\overline{z}^3\overline{\zeta}
+
\tfrac{1}{2}\,z^2\zeta^2\overline{\zeta}^3
-
6\,z^3\zeta\overline{z}^2\overline{\zeta}
-
6\,z^2\zeta\overline{z}^2\overline{\zeta}^2
\\
&
\ \ \ \ \
+
3\,z^2\zeta^2\overline{z}\overline{\zeta}^2
+
2\,z^2\zeta^2\overline{z}^3
+
3\,z\zeta^2\overline{z}^2\overline{\zeta}^2
+
4\,z^3\overline{z}^3\overline{\zeta}
-
2\,z\zeta^2\overline{z}^4,
\endaligned
\]
\[
\aligned
F^8
&
\,=\,
z\overline{z}\zeta^3\overline{\zeta}^3
+
z\zeta^2\overline{z}^2\overline{\zeta}^3
+
z^2\zeta^3\overline{z}\overline{\zeta}^2
+
9\,z^2\zeta^2\overline{z}\overline{\zeta}^3
+
9\,z\zeta^3\overline{z}^2\overline{\zeta}^2
+
\tfrac{14}{3}\,z\zeta^3\overline{z}^3\overline{\zeta}
\\
&
\ \ \ \ \
+
12\,z^2\zeta\overline{z}^4\overline{\zeta}
-
6\,z\zeta^2\overline{z}^4\overline{\zeta}
-
6\,z^4\zeta\overline{z}\overline{\zeta}^2
+
9\,z^2\zeta^2\overline{z}^2\overline{\zeta}^2
+
6\,z^2\zeta\overline{z}^3\overline{\zeta}^2
+
12\,z^4\zeta\overline{z}^2\overline{\zeta}
\\
&
\ \ \ \ \
+
6\,z^3\zeta^2\overline{z}^2\overline{\zeta}
+
\tfrac{14}{3}\,z^3\zeta\overline{z}\overline{\zeta}^3
-
14\,z^3\zeta\overline{z}^2\overline{\zeta}^2
-
14\,z^2\zeta^2\overline{z}^3\overline{\zeta}
-
4\,z^2\zeta\overline{z}^2\overline{\zeta}^3
\\
&
\ \ \ \ \
-\,
4\,z^2\zeta^3\overline{z}^2\overline{\zeta}
-
6\,z^3\zeta^2\overline{z}\overline{\zeta}^2
-
6\,z\zeta^2\overline{z}^3\overline{\zeta}^2
+
\zeta^3\overline{z}^3\overline{\zeta}^2
+
4\,z^3\zeta\overline{z}^4
-
12\,z^4\overline{z}\overline{\zeta}^3
\endaligned
\]
\[
\aligned
{}
&
\ \ \ \ \
+
4\,z^4\overline{z}^3\overline{\zeta}
+
\tfrac{10}{3}\,z^3\overline{z}^2\overline{\zeta}^3
+
3\,z^5\overline{z}^2\overline{\zeta}
+
z^3\zeta\overline{\zeta}^4
+
\zeta^4\overline{z}^3\overline{\zeta}
-
5\,z^4\zeta\overline{z}^3
+
\tfrac{10}{3}\,z^2\zeta^3\overline{z}^3
\\
&
\ \ \ \ \
-\,4\,z^4\overline{z}^2\overline{\zeta}^2
-
5\,z^3\overline{z}^4\overline{\zeta}
-
12\,z\zeta^3\overline{z}^4
-
4\,z^2\zeta^2\overline{z}^4
+
3\,z^2\zeta\overline{z}^5
-
\tfrac{2}{3}\,z^5\overline{\zeta}^3
+
\tfrac{13}{6}\,z^4\overline{\zeta}^4
\\
&
\ \ \ \ \
-\,2\,z^5\overline{z}^3
+
\tfrac{13}{6}\,\zeta^4\overline{z}^4
-
2\,z^3\overline{z}^5
-
\tfrac{2}{3}\,
\zeta^3\overline{z}^5
+
z^3\zeta^2\overline{\zeta}^3.
\endaligned
\]

The general infinitesimal CR automorphism,
depending on $5$ real constants $a, b, c, d, e \in \R$,
is $L = A\, \partial_z + B\, \partial_\zeta + C\, \partial_w$,
where:
\[
\aligned
A^0
&
\,=\,
a
+
\isqrt\,b,
\\
A^1
&
\,=\,
\big(
-a+\isqrt\,b
\big)\,\zeta
+
\big(
2\,a
-
3\,c
+
2\,\isqrt\,b
+
\isqrt\,d
\big)\,z,
\\
A^2
&
\,=\,
\big(
-2\,a
+
2\,c
\big)\,w
+
\big(
2\,a
-
4\,c
+
2\,\isqrt\,d
\big)\,z^2
+
\big(
-c
+
\isqrt\,d
\big)\,z\zeta,
\\
A^3
&
\,=\,
\big(
-2\,a
+
2\,c
\big)\,\zeta w
+
\big(
-4\,a
+
4\,c
\big)\,
zw,
\endaligned
\]
\[
\aligned
A^4
&
\,=\,
0,
\\
A^5
&
\,=\,
0,
\\
A^6
&
\,=\,
0,
\endaligned
\]
where:
\[
\aligned
B^0
&
\,=\,
c
+
\isqrt\,d,
\\
B^1
&
\,=\,
\big(
2\,\isqrt\,d
+
4\,\isqrt\,b
\big)\,\zeta
+
\big(
4\,a
+
4\,\isqrt\,d
\big)\,z,
\\
B^2
&
\,=\,
\big(
8\,a
-
8\,c
+
4\,\isqrt\,d
-
4\,\isqrt\,b
\big)\,z^2
+
\big(
-6\,a
+
6\,\isqrt\,b
-
c
+
\isqrt\,d
\big)\,\zeta^2
\\
&
\ \ \ \ \ 
+
\big(
4\,a
-
12\,c
+
8\,\isqrt\,d
\big)\,z\zeta,
\endaligned
\]
\[
\aligned
B^3
&
\,=\,
\big(
4\,c
-
4\,\isqrt\,d
\big)\,z^3
+
\big(
2\,\isqrt\,b
-
2\,a
\big)\,\zeta^3
+
\big(
-12\,\isqrt\,b
+
12\,a
\big)\,z^2\zeta
+
\big(
-8\,a
+
8\,c
\big)\,\zeta w
\\
&
\ \ \ \ \
+
\big(
-12\,\isqrt\,b
+
12\,a
+
6\,\isqrt\,d
-
6\,c
\big)\,z\zeta^2,
\endaligned
\]
\[
\aligned
B^4
&
\,=\,
\big(
8\,a
-
8\,c
\big)\,
z^2w
+
\big(
-6\,a
+
6\,\isqrt\,b
\big)\,z^4
+
\big(
-12\,\isqrt\,d
+
24\,\isqrt\,b
-
24\,a
+
12\,c
\big)\,z^3\zeta
\\
&
\ \ \ \ \
+
\big(
8\,a
-
2\,c
+
2\,\isqrt\,d
-
8\,\isqrt\,b
\big)\,z\zeta^3
+
\big(
-12\,a
+
12\,\isqrt\,b
+
12\,c
-
12\,\isqrt\,d
\big)\,z^2\zeta^2
\\
&
\ \ \ \ \ 
+
\big(
-12\,a+12\,c
\big)\,\zeta^2w,
\endaligned
\]
\[
\aligned
B^5
&
\,=\,
\big(
-12\,\isqrt\,b
+
12\,a
+
6\,\isqrt\,d
-
6\,c
\big)\,z^5
+
\big(
-24\,c
-
30\,\isqrt\,b
+
30\,a
+
24\,\isqrt\,d
\big)\,z^4\zeta
\\
&
\ \ \ \ \
+
\big(
-20\,a
+
20\,\isqrt\,b
+
8\,c
-
8\,\isqrt\,d
\big)\,z^2\zeta^3
+
\big(
-12\,c
+
12\,\isqrt\,d
\big)\,z^3\zeta^2
+
\big(
-4\,a
+
4\,c
\big)\,\zeta^3w
\\
&
\ \ \ \ \
+
\big(
-24\,c
+
24\,a
\big)\,z\zeta^2w
+
\big(
-24\,c
+
24\,a
\big)\,z^2\zeta w,
\endaligned
\]
and where:
\[
\aligned
C^0
&
\,=\,
\isqrt\,e,
\\
C^1
&
\,=\,
\big(
2\,a
-
2\,\isqrt\,b
\big)\,z,
\\
C^2
&
\,=\,
\big(
4\,a
-
6\,c
\big)\,w
+
\big(
c
-
\isqrt\,d
\big)\,z^2,
\\
C^3
&
\,=\,
\big(
4\,a
-
4\,c
\big)\,zw,
\\
C^4
&
\,=\,
0,
\\
C^5
&
\,=\,
0,
\\
C^6
&
\,=\,
0,
\\
C^7
&
\,=\,
0.
\endaligned
\]
and the related $5$ holomorphic vector fields
$e_1$, $e_2$, $e_3$, $e_4$, $e_5$ have structure:
\[
\aligned
{}
\!\!\!\!\!\!\!\!\!\!\!\!\!\!\!
[e_1,e_2]
&
\,=\,
-\,4\,e_4
-
4\,e_5,
&
\ \ \ \ \
[e_1,e_3]
&
\,=\,
-\,2\,e_1,
&
\ \ \ \ \
[e_1,e_4]
&
\,=\,
2\,e_2
+
4\,e_4,
&
\ \ \ \ \
[e_1,e_5]
&
\,=\,
2\,e_2
-
4\,e_5,
\\
&
&
[e_2,e_3]
&
\,=\,
-\,4\,e_2
-
4\,e_4,
&
\ \ \ \ \
[e_2,e_4]
&
\,=\,
0,
&
\ \ \ \ \
[e_2,e_5]
&
\,=\,
0,
\\
&
&
&
&
\ \ \ \ \
[e_3,e_4]
&
\,=\,
2\,e_4,
&
\ \ \ \ \
[e_3,e_5]
&
\,=\,
-\,2\,e_2
+
6\,e_5,
\\
&
&
&
&
&
&
\ \ \ \ \
[e_4,e_5]
&
\,=\,
0.
\endaligned
\]
\end{Theorem}

This Lie algebra $\mathfrak{g}$ 
has the derived series of dimensions $5$, $4$, $2$, $0$,
with:
\[
[\mathfrak{g},\mathfrak{g}]
\,=\,
\Span\,
\big(
\underline{-4\,e_4-4\,e_5},\,\,\,
-2\,e_1,\,\,\,
\underline{2\,e_2+4\,e_4},\,\,\,
\underline{-4\,e_2-4\,e_4}
\big).
\]
The three underlined vector fields span a $3$-dimensional
Abelian ideal $\mathfrak{a} \subset \mathfrak{g}$,
whose value at the origin $0 \in \C^3$
spans a maximally real $3$-plane.
This is coherent with Fels-Kaup's item {\bf (3)}.

\medskip

Next, assume $F_{3,0,0,2,0} \equiv 0$, or equivalently,
$\frac{1}{4} \overline{\Waux}_0 \equiv 0$. 
Some differential consequences are:
\[
F_{4,0,0,2,0}
\,=\,
0,
\ \ \ \ \ \ \ \ \ \ \ \ \ \ \ \ \ \ \ \
F_{3,0,1,2,0}
\,=\,
0,
\ \ \ \ \ \ \ \ \ \ \ \ \ \ \ \ \ \ \ \
F_{3,0,0,3,0}
\,=\,
0,
\]
hence up to order $6$:
\[
\aligned
u
&
\,=\,
z\overline{z}
+
\tfrac{1}{2}\,
\overline{z}^2\zeta
+
\tfrac{1}{2}\,
z^2\overline{\zeta}
+
z\overline{z}\zeta\overline{\zeta}
+
\tfrac{1}{2}\,
\overline{z}^2\zeta\zeta\overline{\zeta}
+
\tfrac{1}{2}\,
z^2\overline{\zeta}\zeta\overline{\zeta}
+
z\overline{z}\zeta\overline{\zeta}\zeta\overline{\zeta}
\\
&
\ \ \ \ \
+
2\,\Re\,
\Big\{
z^5\overline{\zeta}\,
F_{5,0,0,1,0}
+
z^3\overline{z}^2\overline{\zeta}\,
F_{3,0,2,1,0}
\Big\}
+
{\rm O}_{z,\zeta,\overline{z},\overline{\zeta},v}(7).
\endaligned
\]
Suppose the graphed equation for $M'$ is similar. Then
$F_{5,0,0,1,0}$ is a relative invariant, and it creates a
branching:

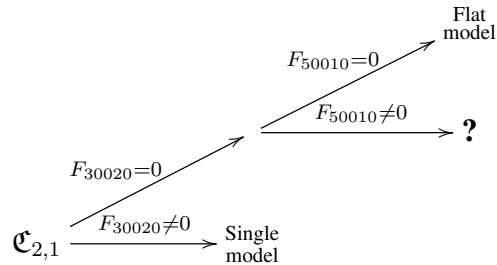
\[
\xymatrix{
&&
&&
{\substack{
\text{Flat}
\\
\text{model}}}
\\
&&
\ar[urr]^{F_{50010}=0}
\ar[rr]^{F_{50010}\neq0}
&&
\text{\bf ?}
\\
\ar[urr]^{F_{30020}=0}
\mathfrak{C}_{2,1}
\ar[rr]^{F_{30020}\neq0}
&&
{\substack{
\text{Single}
\\
\text{model}}}
&&
}
\]
A further sub-branching could be created by the other relative
invariant $F_{3,0,2,1,0}$, but this is not the case. The following
result establishes, by normal forms
techniques, Pocchiola's characterization of the flat model.

\begin{Theorem}
In the branch $F_{3,0,0,2,0} = 0 = F_{5,0,0,1,0}$, 
if $M^5 \in \mathfrak{C}_{2,1}$ is homogeneous, 
then all $G_{h,i,j,k,l} = 0$, 
and $M$ coincides with the Gaussier-Merker
representation of the flat model:
\[
u
\,=\,
\maux
+
0
\,=\,
\frac{z\overline{z}
+\frac{1}{2}\,\overline{z}^2\zeta
+\frac{1}{2}\,z^2\overline{\zeta}}{
1-\zeta\overline{\zeta}}.
\eqno\qed
\]
\end{Theorem}

Thus, in this top-most (degenerate) branch, 
$F_{3,0,2,1,0} = 0$ is {\em implied}, suprisingly.

\medskip

Next, in the branch $F_{3,0,0,2,0} = 0$ and
$F_{5,0,0,1,0} \neq 0$, one can use $\lambda \in \C$
to normalize $F_{5,0,0,1,0} := 1$, so $\lambda = 1$.
The final tree will be explained by the third theorem:
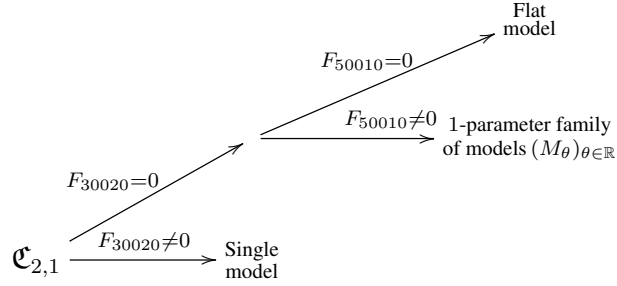
\[
\xymatrix{
&&
&&
{\substack{
\text{Flat}
\\
\text{model}}}
\\
&&
\ar[urr]^{F_{50010}=0}
\ar[rr]^{F_{50010}\neq0}
&&
{\substack{
1\text{-parameter family}
\\
\text{of models}\,(M_\theta)_{\theta\in\R}}}
\\
\ar[urr]^{F_{30020}=0}
\mathfrak{C}_{2,1}
\ar[rr]^{F_{30020}\neq0}
&&
{\substack{
\text{Single}
\\
\text{model}}}
&&
}
\]

\begin{Theorem}
In the branch $F_{3,0,0,2,0} = 0$ and $F_{5,0,0,1,0} = 1$,
three supplementary (real) normalizations hold:
\[
\aligned
F_{6,0,0,1,0}
&
\,:=\,
0,
&
\ \ \ \ \ \ \ \ \ \ \ \ \ \ \ \ \ \ \ \
\text{so}\ \
\alpha
&
\,:=\,
0,
\\
\Im\,F_{4,0,3,0,0}
&
\,:=\,
0,
&
\ \ \ \ \ \ \ \ \ \ \ \ \ \ \ \ \ \ \ \
\text{so}\ \
\rho
&
\,:=\,
0,
\endaligned
\]
so that the isotropy is reduced to be zero-dimensional.
Notably, a constant value
for $F_{3,0,2,1,0} = -\,15$ is also implied.

Furthermore, abbreviating:
\[
\theta
\,:=\,
\Re\,F_{4,0,3,0,0},
\]
which is a {\em free} absolute invariant,
all coefficients $F_{h,i,j,k,l} \in \C$ are
uniquely determined in terms of $\theta \in \R$, with:
\[
F
\,=\,
F^2+F^3+F^4+F^5+F^6+F^7+F^8+F^9+F^{10}
+
{\rm O}_{z,\zeta,\overline{z},\overline{\zeta},v}(11),
\]
where:
\[
\aligned
F^2
&
\,=\,
z\overline{z},
\\
F^3
&
\,=\,
\tfrac{1}{2}\,
\overline{z}^2\zeta
+
\tfrac{1}{2}\,
z^2\overline{\zeta},
\\
F^4
&
\,=\,
z\overline{z}\zeta\overline{\zeta},
\\
F^5
&
\,=\,
\tfrac{1}{2}\,
z^2\zeta\overline{\zeta}^2
+
\tfrac{1}{2}\,
\zeta^2\overline{z}^2\overline{\zeta},
\\
F^6
&
\,=\,
-\,15\,z^3\overline{z}^2\overline{\zeta}
+
z\zeta^2\overline{z}\overline{\zeta}^2
+
\zeta\overline{z}^5
+
z^5\overline{\zeta}
-
15\,z^2\zeta\overline{z}^3,
\endaligned
\]
\[
\aligned
F^7
&
\,=\,
\tfrac{3}{2}\,z^5\overline{\zeta}^2
+
\theta\,z^4\overline{z}^3
-
45\,z^2\zeta\overline{z}^3\overline{\zeta}
+
\tfrac{1}{2}\,z^2\zeta^2\overline{\zeta}^3
+
\theta\,z^3\overline{z}^4
-
\tfrac{15}{2}\,z^4\overline{z}\overline{\zeta}^2
+
5\,z\zeta\overline{z}^4\overline{\zeta}
\\
&
\ \ \ \ \
-\,10\,z^3\overline{z}^2\overline{\zeta}^2
+
\tfrac{3}{2}\,\zeta^2\overline{z}^5
-
45\,z^3\zeta\overline{z}^2\overline{\zeta}
-
\tfrac{15}{2}\,z\zeta^2\overline{z}^4
-
10\,z^2\zeta^2\overline{z}^3
+
\tfrac{1}{2}\,\zeta^3\overline{z}^2\overline{\zeta}^2
\\
&
\ \ \ \ \
+5\,z^4\zeta\overline{z}\overline{\zeta},
\endaligned
\]
\[
\aligned
F^8
&
\,=\,
z\zeta^3\overline{z}\overline{\zeta}^3
-
20\,z^3\zeta^2\overline{z}^2\overline{\zeta}
-
75\,z^3\zeta\overline{z}^2\overline{\zeta}^2
-
75\,z^2\zeta^2\overline{z}^3\overline{\zeta}
-
\tfrac{75}{2}\,z\zeta^2\overline{z}^4\overline{\zeta}
\\
&
\ \ \ \ \
-\,
\tfrac{75}{2}\,z^4\zeta\overline{z}\overline{\zeta}^2
-
20\,z^2\zeta\overline{z}^3\overline{\zeta}^2
-
\tfrac{1}{5}\,\theta\,z^6\overline{z}\overline{\zeta}
+
2\,\theta\,z^4\zeta\overline{z}^3
+
\tfrac{12}{5}\,\theta\,z^5\overline{z}^2\overline{\zeta}
\\
&
\ \ \ \ \
+
3\,\theta\,z^4\overline{z}^3\overline{\zeta}
+
2\,\theta\,z^3\overline{z}^4\overline{\zeta}
+
3\,\theta\,z^3\zeta\overline{z}^4
+
\tfrac{12}{5}\,\theta\,z^2\zeta\overline{z}^5
-
\tfrac{1}{5}\,\theta\,z\zeta\overline{z}^6
-
5\,z^4\overline{z}\overline{\zeta}^3
\\
&
\ \ \ \ \
-\,
5\,z\zeta^3\overline{z}^4
+
5\,\zeta^2\overline{z}^5\overline{\zeta}
+
5\,z^5\zeta\overline{\zeta}^2
-
\tfrac{1}{35}\,\theta\,\zeta\overline{z}^7
-
\tfrac{1}{35}\,\theta\,z^7\overline{\zeta}
-
\tfrac{3}{2}\,z^5\overline{\zeta}^3
-
130\,z^5\overline{z}^3
\\
&
\ \ \ \ \
-\,
\tfrac{325}{6}\,z^4\overline{z}^4
-
130\,z^3\overline{z}^5
-
\tfrac{3}{2}\,\zeta^3\overline{z}^5,
\endaligned
\]
\[
\aligned
F^9
&
\,=\,
\theta\,z^3\overline{z}^4\overline{\zeta}^2
+
\theta\,z^4\zeta^2\overline{z}^3
-
165\,z^2\zeta^2\overline{z}^3\overline{\zeta}^2
-
40\,z^3\zeta\overline{z}^2\overline{\zeta}^3
-
165\,z^3\zeta^2\overline{z}^2\overline{\zeta}^2
\\
&
\ \ \ \ \
-\,
5\,z^4\zeta^2\overline{z}\overline{\zeta}^2
-
5\,z\zeta^2\overline{z}^4\overline{\zeta}^2
-
\tfrac{75}{2}\,z^4\zeta\overline{z}\overline{\zeta}^3
-
\tfrac{75}{2}\,z\zeta^3\overline{z}^4\overline{\zeta}
-
40\,z^2\zeta^3\overline{z}^3\overline{\zeta}
\\
&
\ \ \ \ \
+
\tfrac{18}{5}\,\theta\,z^5\overline{z}^2\overline{\zeta}^2
+
3\,\theta\,z^4\overline{z}^3\overline{\zeta}^2
+
2\,\theta\,z^6\overline{z}\overline{\zeta}^2
+
\tfrac{18}{5}\,\theta\,z^2\zeta^2\overline{\zeta}^5
+
2\,\theta\,z\zeta^2\overline{z}^6
\\
&
\ \ \ \ \
+
3\,\theta\,z^3\zeta^2\overline{\zeta}^4
-
5\,\isqrt\,z^6\overline{\zeta}v
+
5\,\isqrt\,\zeta\overline{z}^6v
+
\tfrac{24}{5}\,\theta\,z^2\zeta\overline{z}^5\overline{\zeta}
+
12\,\theta\,z^4\zeta\overline{z}^3\overline{\zeta}
\\
&
\ \ \ \ \
+
\tfrac{24}{5}\,\theta\,z^5\zeta\overline{z}^2\overline{\zeta}
+
12\,\theta\,z^3\zeta\overline{z}^4\overline{\zeta}
-
\tfrac{1}{5}\,\theta\,z\zeta\overline{z}^6\overline{\zeta}
-
\tfrac{1}{5}\,\theta\,z^6\zeta\overline{z}\overline{\zeta}
-
100\,\isqrt\,z^3\overline{z}^3\overline{\zeta}v
\endaligned
\]
\[
\aligned
{}
&
\ \ \ \ \
-\,30\,\isqrt\,z^5\overline{z}\overline{\zeta}v
-
75\,\isqrt\,z^4\overline{z}^2\overline{\zeta}v
-
\tfrac{335}{3}\,z^5\overline{z}^3\overline{\zeta}
+
5\,z\zeta\overline{z}^7
-
\tfrac{6}{35}\,\theta\,\zeta^2\overline{\zeta}^7
-
\tfrac{335}{3}\,z^3\zeta\overline{z}^5
\\
&
\ \ \ \ \
+
\tfrac{1}{2}\,z^2\zeta^3\overline{\zeta}^4
+
\tfrac{475}{2}\,z^4\overline{z}^4\overline{\zeta}
-
455\,z^2\zeta\overline{z}^6
-
\tfrac{6}{25}\,\theta^2\,z^5\overline{z}^4
-
190\,z^5\zeta\overline{z}^3
-
\tfrac{6}{35}\,\theta\,z^7\overline{\zeta}^2
\\
&
\ \ \ \ \
+
\tfrac{1}{2}\,\zeta^4\overline{z}^2\overline{\zeta}^3
-
455\,z^6\overline{z}^2\overline{\zeta}
+
5\,z^7\overline{z}\overline{\zeta}
-
\tfrac{6}{25}\,\theta^2\,z^4\overline{z}^5
-
190\,z^3\overline{z}^5\overline{\zeta}
-
\tfrac{4}{25}\,\theta^2\,z^6\overline{z}^3
\\
&
\ \ \ \ \
-\,
\tfrac{4}{25}\,\theta^2\,z^3\overline{z}^6
-
\tfrac{15}{2}\,z^5\zeta\overline{\zeta}^3
-
\tfrac{15}{2}\,\zeta^3\overline{z}^5\overline{\zeta}
-
z^5\overline{\zeta}^4
-
\zeta^4\overline{z}^5
+
\tfrac{25}{4}\,\zeta\overline{z}^8
+
\tfrac{25}{4}\,z^8\overline{\zeta}
\\
&
\ \ \ \ \
+
30\,\isqrt\,z\zeta\overline{z}^5v
+
100\,\isqrt\,z^3\zeta\overline{z}^3v
+
75\,\isqrt\,z^2\zeta\overline{z}^4v
+
\tfrac{475}{2}\,z^4\zeta\overline{z}^4,
\endaligned
\]
\[
\aligned
F^{10}
&
\,=\,
z\zeta^4\overline{z}\overline{\zeta}^4
+
\theta\,z^3\zeta^3\overline{z}^4
+
\theta\,z^4\overline{z}^3\overline{\zeta}^3
-
210\,z^3\zeta^2\overline{z}^2\overline{\zeta}^3
-
20\,z^4\zeta\overline{z}\overline{\zeta}^4
+
105\,z^5\zeta\overline{z}^3\overline{\zeta}
\\
&
\ \ \ \ \ 
-\,
210\,z^2\zeta^3\overline{z}^3\overline{\zeta}^2
-
\tfrac{1525}{2}\,z^2\zeta\overline{z}^6\overline{\zeta}
-
\tfrac{1525}{2}\,z^6\zeta\overline{z}^2\overline{\zeta}
-
\tfrac{255}{2}\,z^4\zeta^2\overline{z}\overline{\zeta}^3
\\
&
\ \ \ \ \
+
1725\,z^4\zeta\overline{z}^4\overline{\zeta}
-
20\,z\zeta^4\overline{z}^4\overline{\zeta}
-
70\,z^3\overline{\zeta}^3\overline{z}^2\overline{\zeta}^2
-
\tfrac{255}{2}\,z\zeta^3\overline{z}^4\overline{\zeta}^2
-
70\,z^2\zeta^2\overline{z}^3\overline{\zeta}^3
\\
&
\ \ \ \ \
+
105\,z^3\zeta\overline{z}^5\overline{\zeta}
+
50\,z^7\zeta\overline{z}\overline{\zeta}
+
50\,z\zeta\overline{z}^7\overline{\zeta}
+
\tfrac{3}{175}\,\theta^2\,z\zeta\overline{z}^8
+
\tfrac{3}{175}\,\theta^2\,z^8\overline{z}\overline{\zeta}
\endaligned
\]
\[
\aligned
{}
&
\ \ \ \ \ 
-\,\tfrac{4}{5}\,\theta^2\,z^3\zeta\overline{z}^6
+
2\,\theta\,z\zeta^3\overline{z}^6
+
2\,\theta\,z^6\overline{z}\overline{\zeta}^3
+
\tfrac{12}{5}\,\theta\,z^5\overline{z}^2\overline{\zeta}^3
-
\tfrac{18}{25}\,\theta^2\,z^5\zeta\overline{z}^4
\\
&
\ \ \ \ \
-\,
\tfrac{8}{25}\,\theta^2\,z^6\zeta\overline{z}^3
-
\tfrac{18}{25}\,\theta^2\,z^4\overline{z}^5\overline{\zeta}
-
\tfrac{8}{25}\,\theta^2\,z^3\overline{z}^6\overline{\zeta}
-
\tfrac{72}{175}\,\theta^2\,z^2\zeta\overline{\zeta}^7
-
\tfrac{4}{5}\,\theta^2\,z^6\overline{z}^3\overline{\zeta}
\\
&
\ \ \ \ \
-\,
\tfrac{24}{25}\,\theta^2\,z^5\overline{z}^4\overline{\zeta}
-
\tfrac{24}{25}\,\theta^2\,z^4\zeta\overline{z}^5
+
\tfrac{12}{5}\,\theta\,z^2\zeta^3\overline{z}^5
-
\tfrac{72}{175}\,\theta^2\,z^7\overline{z}^2\overline{\zeta}
-
\tfrac{1}{5}\,\theta\,\zeta^2\overline{z}^7\overline{\zeta}
\\
&
\ \ \ \ \
-\,
\tfrac{1}{5}\,\theta\,z^7\zeta\overline{\zeta}^2
-
20\,\isqrt\,z^6\overline{\zeta}v
+
20\,\isqrt\,\zeta^2\overline{z}^6v
+
24\,\theta\,z^3\zeta^2\overline{z}^4\overline{\zeta}
+
\tfrac{18}{5}\,\theta\,z^6\zeta\overline{z}\overline{\zeta}^2
\endaligned
\]
\[
\aligned
{}
&
\ \ \ \ \
+
\tfrac{12}{5}\,\theta\,z^2\zeta\overline{z}^5\overline{\zeta}^2
+
\tfrac{18}{5}\,\theta\,z\zeta^2\overline{z}^6\overline{\zeta}
+
15\,\theta\,z^4\zeta^2\overline{z}^3\overline{\zeta}
+
15\,\theta\,z^3\zeta\overline{z}^4\overline{\zeta}^2
\\
&
\ \ \ \ \
+
\tfrac{12}{5}\,\theta\,z^5\zeta^2\overline{z}^2\overline{\zeta}
+
18\,\theta\,z^2\zeta^2\overline{z}^5\overline{\zeta}
+
18\,\theta\,z^5\zeta\overline{z}^2\overline{\zeta}^2
+
24\,\theta\,z^4\,\zeta\overline{z}^3\overline{\zeta}^2
\\
&
\ \ \ \ \
-\,
150\,\isqrt\,z^4\overline{z}^2\overline{\zeta}^2v
-
16\,\isqrt\,\theta\,z^3\overline{z}^5v
-
100\,\isqrt\,z^3\overline{z}^3\overline{\zeta}^2v
-
90\,\isqrt\,z^5\overline{z}\overline{\zeta}^2v
+
F_{5,0,5,0,0}\,z^5\overline{z}^5
\\
&
\ \ \ \ \ 
-\,
\tfrac{15}{2}\,z^5\zeta\overline{\zeta}^4
+
5\,\zeta^3\overline{z}^5\overline{\zeta}^2
+
5\,z^5\zeta^2\overline{\zeta}^3
-
150\,z^5\zeta^2\overline{\zeta}^3
+
\tfrac{975}{2}\,z^4\zeta^2\overline{z}^4
+
570\,z^3\zeta^2\overline{z}^5
\endaligned
\]
\[
\aligned
{}
&
\ \ \ \ \
-\,
325\,z^2\zeta^2\overline{z}^6
-
435\,z\zeta^2\overline{z}^7
-
435\,z^7\overline{z}\overline{\zeta}^2
-
325\,z^6\overline{z}^2\overline{\zeta}^2
-
\tfrac{15}{2}\,\zeta^4\overline{z}^5\overline{\zeta}
\\
&
\ \ \ \ \
+
570\,z^5\overline{z}^3\overline{\zeta}^2
+
\tfrac{975}{2}\,z^4\overline{z}^4\overline{\zeta}^2
-
150\,z^3\overline{z}^5\overline{\zeta}^2
+
9\,\theta\,z^4\overline{z}^6
+
\tfrac{164}{7}\,\theta\,z^3\overline{z}^7
+
\tfrac{4}{7}\,\theta\,\zeta^3\overline{z}^7
\\
&
\ \ \ \ \
+
\tfrac{1}{525}\,\theta^2\,\zeta\overline{z}^9
+
\tfrac{4}{7}\,\theta\,z^7\overline{\zeta}^3
+
\tfrac{1}{525}\,\theta^2\,z^9\overline{\zeta}
+
9\,\theta\,z^6\overline{z}^4
+
\tfrac{164}{7}\,\theta\,z^7\overline{z}^3
+
\tfrac{95}{4}\,\zeta^2\overline{z}^8
\\
&
\ \ \ \ \
+
\tfrac{95}{4}\,z^8\overline{\zeta}^2
+
90\,\isqrt\,z\zeta^2\overline{z}^5v
+
100\,\isqrt\,z^3\zeta^2\overline{z}^3v
+
150\,\isqrt\,z^2\zeta^2\overline{z}^4v
+
16\,\isqrt\,\theta\,z^5\overline{z}^3v
\\
&
\ \ \ \ \
-\,
30\,\isqrt\,z^5\zeta\overline{z}\overline{\zeta}v
-
150\,\isqrt\,z^2\zeta\overline{z}^4\overline{\zeta}v
+
150\,\isqrt\,z^4\zeta\overline{z}^2\overline{\zeta}v
+
30\,\isqrt\,z\zeta\overline{z}^5\overline{\zeta}v.
\endaligned
\]

The general infinitesimal CR automorphism,
depending on $5$ real constants $a, b, c, d, e \in \R$,
is $L = A\, \partial_z + B\, \partial_\zeta + C\, \partial_w$,
where:
\[
\aligned
A^0
&
\,=\,
a+\isqrt\,b,
\\
A^1
&
\,=\,
\big(-c+\isqrt\,d\big)\,z
+
\big(-a+\isqrt\,b\big)\,\zeta,
\\
A^2
&
\,=\,
\big(
\tfrac{2}{5}\,\theta\,a
+
5\,\isqrt\,e
\big)\,z^2
+
\big(
-\tfrac{2}{5}\,\theta\,a
+
5\,\isqrt\,e
\big)\,w
+
\big(
-c
+
\isqrt\,d
\big)\,z\zeta,
\\
A^3
&
\,=\,
\big(
-10\,a
-
10\,\isqrt\,b
\big)\,z^3
+
\big(
10\,\isqrt\,b
+
30\,a
\big)\,zw
+
\big(
-\tfrac{2}{5}\,\theta\,a
-
5\,\isqrt\,e
\big)\,\zeta w,
\endaligned
\]
\[
\aligned
A^4
&
\,=\,
\big(
-10\,c
-
5\,\isqrt\,d
\big)\,w^2
+
\big(
10\,a
+
10\,\isqrt\,b
\big)\,z\zeta w
+
\big(
-20\,c
+
10\,\isqrt\,d
\big)\,z^2w,
\\
A^5
&
\,=\,
\big(
-\tfrac{4}{5}\,\theta\,a
-
10\,\isqrt\,e
\big)\,z^5
+
\big(
-5\,c
+
5\,\isqrt\,d
\big)\,z^4\zeta
+
\big(
4\,\theta\,a
-
50\,\isqrt\,e
\big)\,
z^3w
\\
&
\ \ \ \ \
+
\big(
75\,\isqrt\,e
-
6\,\theta\,a
\big)\,zw^2
+
\big(
10\,c
-
5\,\isqrt\,d
\big)\,\zeta w^2,
\endaligned
\]
\[
\aligned
A^6
&
\,=\,
\big(
-20\,a
-
20\,\isqrt\,b
-
\tfrac{1}{5}\,\isqrt\,\theta\,d
+
\tfrac{1}{5}\,\theta\,c
\big)\,z^6
+
\big(
-\tfrac{200}{3}\,a
+
\tfrac{100}{3}\,\isqrt\,b
\big)\,w^3
\\
&
\ \ \ \ \
+
\big(
-2\,\theta\,a
+
25\,\isqrt\,e
\big)\,z\zeta w^2
+
\big(
200\,a
+
100\,\isqrt\,b
\big)\,z^2w^2,
\endaligned
\]
\[
\aligned
A^7
&
\,=\,
\big(
10\,\isqrt\,d
-
10\,c
+
\tfrac{2}{7}\,\isqrt\,\theta\,e
+
\tfrac{4}{175}\,\theta^2\,a
\big)\,z^7
+
\big(
100\,c
+
50\,\isqrt\,d
\big)\,z^3w^2
\\
&
\ \ \ \ \ 
+
\big(
-\tfrac{200}{3}\,a
-
\tfrac{100}{3}\,\isqrt\,b
\big)\,
\zeta w^3
+
\big(
-\tfrac{1}{5}\,\isqrt\,\theta\,d
+
\tfrac{1}{5}\,\theta\,c
\big)\,z^6\zeta
+
\big(
50\,a
+
50\,\isqrt\,b
\big)\,z^4\zeta w
\\
&
\ \ \ \ \
+
\big(
-50\,\isqrt\,d
+
70\,c
\big)\,
z^5w
+
\big(
-100\,c
-
50\,\isqrt\,d
+
20\,\isqrt\,\theta\,e
\big)\,
zw^3,
\endaligned
\]
\[
\aligned
A^8
&
\,=\,
A_{0,0,4}\,w^4
+
\big(
-\tfrac{3}{175}\,\theta^2\,c
+
\tfrac{4}{7}\,\isqrt\,\theta\,b
-
\tfrac{31}{7}\,\theta\,a
+
\tfrac{3}{175}\,\isqrt\,\theta^2\,d
-
\tfrac{125}{2}\,\isqrt\,e
\big)\,z^8
\\
&
\ \ \ \ \
+
\big(
8\,\overline{A_{0,0,4}}
+
\tfrac{1}{2}\,
\overline{B_{1,0,3}}
\big)\,
z^2w^3
+
\big(
-30\,\isqrt\,d
+
30\,c
\big)\,z^5\zeta w
+
\big(
-\tfrac{100}{3}\,c
-
\tfrac{50}{3}\,\isqrt\,d
\big)\,z\zeta w^3
\\
&
\ \ \ \ \
+
\big(
-50\,c
+
50\,\isqrt\,d
\big)\,z^7\zeta
+
\big(
-2\,\isqrt\,\theta\,b
+
10\,\theta\,a
-
50\,\isqrt\,e
\big)\,
z^6w,
\endaligned
\]
where:
\[
\aligned
B^0
&
\,=\,
c+\isqrt\,d,
\\
B^1
&
\,=\,
\big(
\tfrac{4}{5}\,\theta\,a
-
10\,\isqrt\,e
\big)\,z
+
\big(
2\,\isqrt\,d
\big)\,\zeta,
\\
B^2
&
\,=\,
\big(
-40\,\isqrt\,b
-
60\,a
\big)\,z^2
+
\big(
-c
+
\isqrt\,d
\big)\,\zeta^2
+
\big(
10\,a
-
10\,\isqrt\,b
\big)\,w
+
\big(
\tfrac{4}{5}\,\theta\,a
+
10\,\isqrt\,e
\big)\,
z\zeta,
\\
B^3
&
\,=\,
\big(
-30\,\isqrt\,d
+
30\,c
\big)\,z^3
+
\big(
40\,c
+
20\,\isqrt\,d
\big)\,zw
+
\big(
60\,a
\big)\,\zeta w
+
\big(
40\,a
-
140\,\isqrt\,b
\big)\,
z^2\zeta,
\endaligned
\]
\[
\aligned
B^4
&
\,=\,
\big(
-14\,\theta\,a
+
100\,\isqrt\,e
+
6\,\isqrt\,\theta\,b
\big)\,z^4
+
\big(
2\,\theta\,a
+
25\,\isqrt\,e
\big)\,w^2
+
\big(
-40\,c
+
20\,\isqrt\,d
\big)\,z\zeta w
\\
&
\ \ \ \ \
+
\big(
24\,\theta\,a
-
300\,\isqrt\,e
\big)\,z^2w
+
\big(
10\,a
+
10\,\isqrt\,b
\big)\,\zeta^2w
+
\big(
-90\,\isqrt\,d
+
90\,c
\big)\,z^3\zeta
\\
&
\ \ \ \ \ 
+
\big(
-60\,\isqrt\,b
+
60\,a
\big)\,z^2\zeta^2,
\endaligned
\]
\[
\aligned
B^5
&
\,=\,
\big(
-860\,\isqrt\,b
+
900\,a
+
\tfrac{24}{5}\,\isqrt\,\theta\,d
-
\tfrac{24}{5}\,\theta\,c
\big)\,
z^5
+
\big(
40\,c
-
40\,\isqrt\,d
\big)\,
z^3\zeta^2
\\
&
\ \ \ \ \
+
\big(
12\,\isqrt\,\theta\,b
-
20\,\theta\,a
-
100\,\isqrt\,e
\big)\,
z^4\zeta
+
\big(
-300\,a
-
300\,\isqrt\,b
\big)\,z^3w
\\
&
\ \ \ \ \
+
\big(
400\,a
-
200\,\isqrt\,b
\big)\,zw^2
+
\big(
150\,\isqrt\,e
\big)\,\zeta w^2
+
\big(
56\,\theta\,a
+
200\,\isqrt\,e
\big)\,z^2\zeta w,
\endaligned
\]
\[
\aligned
B^6
&
\,=\,
\big(
-770\,\isqrt\,d
+
690\,c
+
4\,\isqrt\,\theta\,e
+
\tfrac{32}{25}\,\theta^2\,a
-
\tfrac{24}{25}\,\isqrt\,\theta^2\,b
\big)\,
z^6
\\
&
\ \ \ \ \
+
\big(
-\tfrac{100}{3}\,c
+
\tfrac{50}{3}\,\isqrt\,d
\big)\,w^3
+
\big(
-250\,c
-
350\,\isqrt\,d
-
30\,\isqrt\,\theta\,e
-
\tfrac{12}{5}\,\theta^2\,a
\big)\,z^4w
\\
&
\ \ \ \ \
+
\big(
400\,a
+
200\,\isqrt\,b
\big)\,
z\zeta w^2
+
\big(
-2\,\theta\,a
+
25\,\isqrt\,e
\big)\,
\zeta^2 w^2
+
\big(
-6\,\theta\,a
+
6\,\isqrt\,\theta\,b
\big)\,
z^4\zeta^2
\\
&
\ \ \ \ \
+
\big(
600\,c
+
300\,\isqrt\,d
-
60\,\isqrt\,\theta\,e
\big)\,z^2w^2
+
\big(
-1500\,a
-
300\,\isqrt\,b
\big)\,z^3\zeta w
\\
&
\ \ \ \ \
+
\big(
300\,\isqrt\,e
+
24\,\theta\,a
\big)\,
z^2\zeta^2w
+
\big(
\tfrac{48}{5}\,\isqrt\,\theta\,d
-
1360\,\isqrt\,b
+
920\,a
-
\tfrac{48}{5}\,\theta\,c
\big)\,z^5\zeta,
\endaligned
\]
\[
\aligned
B^7
&
\,=\,
\big(
\tfrac{1056}{7}\,\isqrt\,\theta\,b
-
168\,\theta\,a
+
\tfrac{144}{175}\,\theta^2\,c
-
\tfrac{144}{175}\,\isqrt\,\theta^2\,d
\big)\,z^7
\\
&
\ \ \ \ \
+
\big(
-48\,\overline{A_{0,0,4}}
-
6\,\overline{B_{1,0,3}}
+
60\,\theta\,a
-
750\,\isqrt\,e
\big)\,z^3w^2
+
\big(
-200\,c
\big)\,\zeta w^3
\\
&
\ \ \ \ \
+
\big(
\tfrac{56}{25}\,\theta^2\,a
-
\tfrac{48}{25}\,\isqrt\,\theta^2\,b
+
1460\,c
-
1460\,\isqrt\,d
+
4\,\isqrt\,\theta\,e
\big)\,z^6\zeta
\\
&
\ \ \ \ \
+
\big(
-\tfrac{24}{5}\,\theta^2\,a
-
60\,\isqrt\,\theta\,e
+
100\,c
+
100\,\isqrt\,d
\big)\,z^4\zeta w
\\
&
\ \ \ \ \
+
\big(
\tfrac{24}{25}\,\isqrt\,\theta\,d
-
\tfrac{24}{25}\,\theta\,c
+
900\,a
-
900\,\isqrt\,b
\big)\,z^5\zeta^2
\\
&
\ \ \ \ \
+
\big(
360\,\theta\,a
+
144\,\isqrt\,\theta\,b
+
5100\,\isqrt\,e
\big)\,z^5w
+
\big(
B_{1,0,3}
\big)\,zw^3
\\
&
\ \ \ \ \
+
\big(
-1000\,a
+
200\,\isqrt\,b
\big)\,z^3\zeta^2 w
+
\big(
-400\,c
+
700\,\isqrt\,d
\big)\,z^2\zeta w^2,
\endaligned
\]
and where:
\[
\aligned
C^0
&
\,=\,
\isqrt\,e,
\\
C^1
&
\,=\,
\big(
2\,a
-
2\,\isqrt\,b
\big)\,z,
\\
C^2
&
\,=\,
\big(
c
-
\isqrt\,d
\big)\,z^2
+
\big(
-2\,c
\big)\,w,
\\
C^3
&
\,=\,
\big(
\tfrac{4}{5}\,\theta\,a
+
10\,\isqrt\,e
\big)\,zw,
\endaligned
\]
\[
\aligned
C^4
&
\,=\,
\big(
10\,\isqrt\,b
\big)\,w^2
+
\big(
-10\,a
-
10\,\isqrt\,b
\big)\,z^2w,
\\
C^5
&
\,=\,
\big(
2\,c
-
2\,\isqrt\,d
\big)\,
z^5
+
\big(
-20\,c
+
10\,\isqrt\,d
\big)\,zw^2,
\\
C^6
&
\,=\,
\big(
-4\,\theta\,a
\big)\,
w^3
+
\big(
2\,\theta\,a
-
25\,\isqrt\,e
\big)\,z^2w^2,
\endaligned
\]
\[
\aligned
C^7
&
\,=\,
\big(
\tfrac{2}{35}\,\isqrt\,\theta\,d
-
\tfrac{2}{35}\,\theta\,c
\big)\,z^7
+
\big(
-20\,a
-
20\,\isqrt\,b
\big)\,z^5w
\\
&
\ \ \ \ \
+
\big(
\tfrac{400}{3}\,a
+
\tfrac{200}{3}\,\isqrt\,b
\big)\,zw^3,
\endaligned
\]
\[
\aligned
C^8
&
\,=\,
\big(
-25\,\isqrt\,d
+
10\,\isqrt\,\theta\,e
\big)\,w^4
+
\big(
-\tfrac{25}{2}\,\isqrt\,d
+
\tfrac{25}{2}\,c
\big)\,z^8
+
\big(
\tfrac{100}{3}\,c
+
\tfrac{50}{3}\,\isqrt\,d
\big)\,z^2w^3
\\
&
\ \ \ \ \
+
\big(
10\,\isqrt\,d
-
10\,c
\big)\,z^6w,
\endaligned
\]
\[
\aligned
C^9
&
\,=\,
\big(
-\tfrac{2}{525}\,\isqrt\,\theta^2\,d
+
\tfrac{2}{525}\,\theta^2\,c
\big)\,z^9
+
\big(
\tfrac{4}{7}\,\isqrt\,\theta\,b
+
\tfrac{4}{7}\,\theta\,a
\big)\,z^7w
\\
&
\ \ \ \ \
+
\big(
4\,\theta\,a
-
50\,\isqrt\,e
\big)\,z^5w^2
+
\big(
2\,\overline{A_{0,0,4}}
\big)\,zw^4,
\endaligned
\]
and the related $5$ holomorphic vector fields
$e_1$, $e_2$, $e_3$, $e_4$, $e_5$ have structure:
\[
\aligned
{}
\!\!\!\!\!\!\!\!\!\!\!\!\!\!\!
[e_1,e_2]
&
\,=\,
-\,\tfrac{4}{5}\,\theta\,
e_4
-
4\,e_5,
&
\ \ \ \ \
[e_1,e_3]
&
\,=\,
0,
&
\ \ \ \ \
[e_1,e_4]
&
\,=\,
2\,e_2,
&
\ \ \ \ \
[e_1,e_5]
&
\,=\,
\tfrac{2}{5}\,\theta\,\,e_2
-
20\,e_4,
\\
&
&
[e_2,e_3]
&
\,=\,
-\,2\,e_2,
&
\ \ \ \ \
[e_2,e_4]
&
\,=\,
0,
&
\ \ \ \ \
[e_2,e_5]
&
\,=\,
0,
\\
&
&
&
&
\ \ \ \ \
[e_3,e_4]
&
\,=\,
2\,e_4,
&
\ \ \ \ \
[e_3,e_5]
&
\,=\,
2\,e_5,
\\
&
&
&
&
&
&
\ \ \ \ \
[e_4,e_5]
&
\,=\,
0.
\endaligned
\]
\end{Theorem}

This Lie algebra $\mathfrak{g}$ 
has the derived series of dimensions $5$, $3$, $0$,
with:
\[
[\mathfrak{g},\mathfrak{g}]
\,=\,
\Span\,
\big(
-\tfrac{4}{5}\,\theta\,e_4-4\,e_5,\,\,\,
2\,e_2,\,\,\,
\tfrac{2}{5}\,\theta\,e_2-20\,e_4
\big).
\]
These three vector fields form a $3$-dimensional
Abelian ideal $\mathfrak{a} \subset \mathfrak{g}$,
whose value at the origin $0 \in \C^3$
spans a maximally real $3$-plane.
This is coherent with Fels-Kaup's items {\bf (2a)},
{\bf (2b)}, {\bf (2c)}.

\medskip

The proofs rely upon studying the {\sl fundamental equation:}
\[
\aligned
0
&
\,=\,
-\,u'
+
F'\big(z',\zeta',\overline{z}',\overline{\zeta}',v'\big),
\endaligned
\]
where:
\[
\aligned
z'
&
\,=\,
f(z,\zeta,w),
\\
\zeta'
&
\,=\,
g(z,\zeta,w),
\\
w'
&
\,=\,
h(z,\zeta,w),
\endaligned
\]
for $(z, \zeta, w) \in M$, namely for $u = F(z, \zeta,
\overline{z}, \overline{\zeta}, v)$.

\medskip

These power series normal forms confirm Fels-Kaup's
main result that there is a one-to-one correspondence
between affine equivalence classes of
homogeneous parabolic surfaces $S^2 \subset \R^3$,
and biholomorphic equivalence classes 
of CR homogeneous $\mathfrak{C}_{2,1}$ hypersurfaces $M^5 \subset \C^3$.
In the thickest branch $F_{3,0,0,2,0} = 0 \neq F_{5,0,0,1,0}$, 
homogeneous models both depend upon a (free) real 
parameter $\theta \in \R$. Therefore, 
in this branch,
the tubes {\bf (2a)}, {\bf (2b)}, {\bf (2c)}
are closed representations 
of models (a computation
of $F_{3,0,0,2,0}$ and of $F_{5,0,0,1,0}$ confirms this).
From the implicit equations {\bf (2a)}, {\bf (2b)}, {\bf (2c)}, some
(non-closed) graphed power series can be written.

\medskip

A detailed presentation of all (very lengthy) manuscript 
calculations would 
extend the size of this article beyond 50 further pages,
{\em cf.}~{\cite{Foo-Merker-2019, Foo-Merker-Ta-2020, Merker-2021}}.
We would be very interested in
knowing whether other authors can confirm these results using
alternative methods.
By the way, the coefficients
$A_{0,0,4}$, $B_{1,0,3}$, $F_{5,0,5,0,0}$ in
the third theorem can be determined [this requires to push computations
up to order $11$].

\medskip\noindent{\bf Acknowledgments.}
Dennis The provided quick help about the existence of Abelian ideals. 



\vfill
\begin{thebibliography}{XL}

{\scriptsize

{\bf\bibitem{Abdalla-Dillen-Vrancken-1997}
{\rm Abdalla}}, B.; {\rm Dillen}, F.; {\rm Vrancken}, L.:
{\em Affine homogeneous surfaces in $\R^3$ 
with vanishing Pick invariant},
Abh. Math. Sem. Univ. Hamburg {\bf 67} (1997), 105--115.

\smallskip

{\bf\bibitem{Chen-Foo-Merker-Ta-2020}
{\rm Chen}}, Z.; {\rm Foo}, W.G.; {\rm Merker}, J.; {\rm Ta}, T.A.:
{\em Normal forms for rigid $\mathfrak{C}_{2,1}$ hypersurfaces 
$M^5 \subset \C^3$},
to appear in Taiwanese Journal of Mathematics,
{\tiny\sf arxiv.org/abs/1912.01655/}, 32 pages.

\smallskip

{\bf\bibitem{Chen-Merker-2019}
{\rm Chen}}, Z.; {\rm Merker}, J.:
{\em On differential invariants of parabolic surfaces},
{\tiny\sf arxiv.org/abs/1908.07867/},
to appear in Dissertationes Mathematic{\ae} 2021,
doi:10.4064/dm816-8-2020, 110 pages.

\smallskip

{\bf\bibitem{Doubrov-Komrakov-Rabinovich-1996}
{\rm Doubrov}}, B.; Komrakov, B.; Rabinovich, M.: 
{\em Homogeneous surfaces in the three-dimensional affine geometry}. 
Geometry and topology of submanifolds, 
VIII (Brussels, 1995/Nordfjordeid, 1995), 
168--178, World Sci. Publ., River Edge, NJ, 1996.

\smallskip

{\bf\bibitem{Doubrov-Merker-The-2020}
{\rm Doubrov}}, B.; Merker, J.; The, D.:
{\em Classification of simply-transitive Levi non-degenerate
hypersurfaces in $\C^3$}, 
{\tiny\sf arxiv.org/abs/2010.06334/}, 31 pages.

\smallskip

{\bf\bibitem{Eastwood-Ezhov-1999}
{\rm Eastwood}}, M.; {\rm Ezhov}, V.:
{\em On affine normal forms and a classification 
of homogeneous surfaces in affine three-space}, 
Geom. Dedicata {\bf 77} (1999), no.~1, 11--69.

\smallskip

{\bf\bibitem{Fels-Kaup-2007}
{\rm Fels}}, M.; {\rm Kaup}, W.:
{\em CR manifolds of dimension $5$: a Lie algebra approach},
J. Reine Angew. Math. {\bf 604} (2007), 47--71.

\smallskip

{\bf\bibitem{Fels-Kaup-2008}
{\rm Fels}}, M.; {\rm Kaup}, W.:
{\em Classification of Levi degenerate homogeneous CR-manifolds in 
dimension $5$},
Acta Math. {\bf 201} (2008), 1--82.

\smallskip

{\bf\bibitem{Foo-Merker-2019}
{\rm Foo}}, W.G.; {\rm Merker}, J.:
{\em Differential $\{e\}$-structures for equivalences
of $2$-nondegenerate Levi rank $1$ hypersurfaces $M^5 \subset \C^3$},
{\tiny\sf arxiv.org/abs/1901.02028/}, 72 pages.

\smallskip

{\bf\bibitem{Foo-Merker-Ta-2019}
{\rm Foo}}, W.G.; {\rm Merker}, J.; {\rm Ta}, T.-A.:
{\em Rigid equivalences of $5$-dimensional $2$-nondegenerate
rigid real hypersurfaces $M^{5}\subset\mathbb{C}^{3}$ 
of constant Levi rank $1$},
{\tiny\sf arxiv.org/abs/1904.02562/},
32 pages, Michigan Math. J., to appear.

\smallskip

{\bf\bibitem{Foo-Merker-Ta-2020}
{\rm Foo}}, W.G.; {\rm Merker}, J.; {\rm Ta}, T.-A.:
{\em On convergent Poincar\'e-Moser reduction
for Levi degenerate embedded $5$-dimensional CR manifolds},
{\tiny\sf arxiv.org/abs/2003.01952/},
71 pages.

\smallskip

{\bf\bibitem{Gaussier-Merker-2003}
{\rm Gaussier}}, H.; {\rm Merker}, J.:
{\em A new example of uniformly Levi degenerate hypersurface in 
$\C^3$}, Ark. Mat. {\bf 41} (2003), no.~1, 85--94.
Erratum: {\bf 45} (2007), no.~2, 269--271.

\smallskip

{\bf\bibitem{Isaev-2011}
{\rm Isaev}}, A.:
{\em Spherical tube hypersurfaces},
Lecture Notes in Mathematics, 2020, 
Springer, Heidelberg, 2011, xii+220~pp.

\smallskip

{\bf\bibitem{Isaev-2016}
{\rm Isaev}}, A.:
{\em Affine rigidity of Levi degenerate tube hypersurfaces},
J. Differential Geom. {\bf 104} (2016), no.~1, 111--141.

\smallskip

{\bf\bibitem{Isaev-2016-bis}
{\rm Isaev}}, A.:
{\em On the CR-curvature of Levi degenerate tube hypersurfaces},
{\tiny\sf arxiv.org/abs/1608.02919/}, 2016, 11 pages.

\smallskip

{\bf\bibitem{Isaev-2018}
{\rm Isaev}}, A.:
{\em Zero CR-curvature equations for Levi degenerate hypersurfaces
via Pocchiola's invariants}, 
Ann. Fac. Sci. Toulouse {\bf 28} (2019) no.~5, 957--976.

\smallskip

{\bf\bibitem{Kolar-Kossovskiy-2019}
{\rm Kolar}}, M.; {\rm Kossovskiy}, I.:
{\em A complete normal form for everywhere Levi
degenerate hypersurfaces in $\C^3$}, 
{\tiny\sf arxiv.org/abs/1905.05629/}

\smallskip

{\bf\bibitem{Jacobowitz-1990}
{\rm Jacobowitz}}, H.:
{\em An introduction to CR structures}, 
Math. Surveys and Monographs, 32, Amer. Math. Soc., Providence, 1990,
x+237~pp.

\smallskip

{\bf\bibitem{Loboda-2020} 
{\rm Loboda}}, A.V.:
{\em Holomorphically Homogeneous Real Hypersurfaces in $\C^3$}
(Russian),
to appear in the Proceedings of the Moscow Mathematical Society,
{\tiny\sf arxiv.org/abs/2006.07835/}, 2020, 56 pages.

\smallskip

{\bf\bibitem{Medori-Spiro-2014}
{\rm Medori}}, C.; {\rm Spiro}, A.: 
{\em The equivalence problem for 5-dimensional Levi degenerate 
CR manifolds}, Int. Math. Res. Not. IMRN {\bf 2014},
no.~20, 5602--5647. 

\smallskip

{\bf\bibitem{Merker-2021}
{\rm Merker}}, J.:
{\em Equivalences of PDE systems associated to 
degenerate para-CR structures: foundational aspects},
{\tiny\sf arxiv.org/abs/2101.05559/}

\smallskip

{\bf\bibitem{Merker-Nurowski-2019}
{\rm Merker}}, J.; {\rm Nurowski}, P.:
{\em On degenerate para-CR structures: 
Cartan reduction and homogeneous models},
{\tiny\sf arxiv.org/abs/2003.08166/} (2020), 37 pages.

\smallskip

{\bf\bibitem{Merker-Nurowski-2021}
{\rm Merker}}, J.; {\rm Nurowski}, P.:
{\em Homogeneous CR and para-CR structures in dimensions $5$ and $3$},
to appear.

\smallskip

{\bf\bibitem{Merker-Pocchiola-2018}
{\rm Merker}}, J.; {\rm Pocchiola}, S.:
{\em Explicit absolute parallelism for $2$-nondegenerate real
hypersurfaces $M^5 \subset \C^3$ of constant Levi rank $1$},
Journal of Geometric Analysis, {\bf 30} (2020),
2689--2730, {\tiny\sf 10.1007/s12220-018-9988-3}.
Addendum: 3233--3242,
{\tiny\sf 10.1007/s12220-019-00195-2}.

\smallskip

{\bf\bibitem{Nurowski-Sparling-2003} 
{\rm Nurowski}}, P.; {\rm Sparling}, G.:
{\em Three-dimensional Cauchy-Riemann structures and
second order ordinary differential equations},
Classical Quantum Gravity {\bf 20} (2003), no.~23, 4995--5016.

\smallskip

{\bf\bibitem{Nurowski-Tafel-1988}
{\rm Nurowski}}, P.; {\rm Tafel}, J.:
{\em Symmetries of Cauchy-Riemann spaces},
Lett. Math. Phys. {\bf 15} (1988), no.~1, 31--38. 

\smallskip

{\bf\bibitem{Olver-2018}
{\rm Olver}}, P.J.:
{\em Normal forms for submanifolds under group actions},
Symmetries, differential equations and applications, 1--25.
Springer Proc. Math. Stat. {\bf 266}, Springer, Cham, 2018.

}

\end{thebibliography}
\end{document}